# THICK-THIN DECOMPOSITION FOR QUADRATIC DIFFERENTIALS

KASRA RAFI

## 1. INTRODUCTION

Quadratic differentials arise naturally in the study of Teichmüller space and Teichmüller geodesics, in particular. A quadratic differential $q$ on a surface $S$ defines a singular Euclidean metric on $S$. The explicit nature of singular Euclidean metrics often makes estimates simple and combinatorial arguments possible ([GM91], [Mas93]). On the other hand, by the uniformization theorem, there also exists a canonical hyperbolic metric $\sigma$ in the conformal class of $q$ and the problem of understanding the relationship between the two metrics often arises (see for example [Min92], [Raf05b] and [Raf05a]).

There is a well-known decomposition of a hyperbolic surface $S$ into thick and thin parts ([Thu86], [BP92]). The components of the thin part have a simple topology; they are homeomorphic to annuli. The components of the thick part, on the other hand, have bounded geometry, i.e., the diameter and the injectivity radius of a thick piece are bounded both above and below by constants depending on the topology of $S$ only. The geometry of a thick piece $Y$ is coarsely determined by the topology of a short marking (see [MM00] and Section 3.1 below) in $Y$. For example, the hyperbolic length of a curve $\alpha$ is comparable (up to multiplicative constants depending on the topology of $S$) to the geometric intersection number between $\alpha$ and a given short marking ([Min93], see also (3)).

We are interested in comparing these two metrics restricted to a thick piece $Y$ of $S$. In general the diameter and the area of $Y$ may be much smaller in the quadratic differential metric compared with the hyperbolic metric (see the example at the end of this paper). However, we show that, after scaling appropriately, $Y$ equipped with the quadratic differential metric has a geometry comparable with the geometry of $Y$ equipped with the hyperbolic metric. The following is our main theorem:

**Theorem 1.** *For every thick piece $Y$ of $S$, there exists a constant $\lambda(Y)$ such that, for every essential curve $\alpha$ in $Y$, the $q$–length of $\alpha$ is equal to $\lambda(Y)$ times the $\sigma$–length of $\alpha$ up to multiplicative constants depending only on the topology of $S$.*

This theorem provides a concrete description of the structure of a quadratic differential metric; the surface is decomposed into thick pieces that are glued along flat annuli (see Section 2). In addition, our results have applications to Teichmüller theory. In particular, they have been applied in [CRS06] to study lines of minima in comparison with their corresponding Teichmüller geodesics.

The paper is organized as follows. In Section 2 we briefly review some facts about quadratic differentials that we shall need. In Section 3 we define the *size* $\lambda(Y)$ of a subsurface $Y$, which provides the scale factor in the main theorem. We show in Lemma 3 that a curve whose $\sigma$–length is bounded has $q$–length comparable to $\lambda(Y)$ and prove Theorem 4, which states that the $q$–diameter of $Y$ is comparable to $\lambda(Y)$ and further that the $q$–area of $Y$ is bounded by $\lambda(Y)^2$. The results of Section 3 are used in Section 4 to prove the main theorem. In the last section, we construct an example illustrating Theorem 4 in which the $q$–area is zero while the diameter is comparable to $\lambda(Y)$.

1.1. **Acknowledgment.** I would like to thank Young-Eun Choi and Caroline Series for their careful reading of this manuscript and many helpful suggestions.

## 2. Quadratic differentials

Let $q = \phi(z)\,dz^2$ be a meromorphic quadratic differential of area one on $S$. (For definition and details, see [Str80] and [GL00].) We call the singular Euclidean metric $|q|$ the *q–metric* on $S$. We assume that $q$ has a discrete set of finite critical points (i.e., critical points of $q$ are either zeroes or poles of order 1) and that the surface is punctured at the poles. Every curve $\alpha$ in $S$ can be represented by a geodesic in this metric. This representative might "pass through" the poles even though the poles are removed from the surface. Following the discussion in [Raf05b], we can ignore this difficulty and treat these special geodesics as we would any other geodesic.

2.1. **Length of a curve.** By a *curve* we always mean a non-trivial non-peripheral piecewise-smooth simple closed curve. For a curve $\alpha$ in $S$, the $q$–geodesic representative of $\alpha$ is unique except for the case where it is one of the continuous family of closed geodesics in a flat annulus, which we refer to as the flat annulus corresponding to $\alpha$ (see below). We denote the $q$–length of $\alpha$ by $l_q(\alpha)$ and the $q$–length of the $q$–geodesic representative of $\alpha$ by $l_q([\alpha])$. In general, for any metric $\tau$, $l_\tau(\alpha)$ represents the $\tau$–length of $\alpha$ and $l_\tau([\alpha])$ represents the $\tau$–length of the $\tau$-geodesic representative of $\alpha$.

2.2. **Curvature of a boundary curve of a subsurface.** Let $Y$ be a subsurface of $S$ and $\gamma$ be a boundary component of $Y$. The curvature in $q$ of $\gamma$ with respect to $Y$, $\kappa_Y(\gamma)$, is well defined as a measure with atoms at the corners. We choose the sign to be positive when the acceleration vector points into $Y$. If $\gamma$ is curved non-negatively (or non-positively) with respect to $Y$ at every point, we say it is *monotonically curved* with respect to $Y$. Let $A$ be an annulus in $S$ with boundaries $\gamma_0$ and $\gamma_1$. Suppose both boundaries are monotonically curved with respect to $A$ and $\kappa_A(\gamma_0) \leq 0$. Further, suppose that the boundaries are equidistant from each other, and the interior of $A$ contains no zeroes. We call $A$ a *primitive* annulus and write $\kappa(A) = -\kappa_A(\gamma_0)$. When $\kappa(A) = 0$, $A$ is called a *flat* annulus and is foliated by closed Euclidean geodesics homotopic to the boundaries. See [Min92] for more details.

2.3. **Subsurfaces with geodesic boundary.** Let $\Gamma$ be the set of boundary curves of a subsurface $Y$. For $\gamma \in \Gamma$, let $F_\gamma$ be the flat annulus containing all $q$–geodesic representatives of $\gamma$. Note that $F_\gamma$ is degenerate when the geodesic representative of $\gamma$ is unique. There exists a unique subsurface $\mathsf{Y}$ in the homotopy class of $Y$ that has $q$–geodesic boundaries and is disjoint from the interior of $F_\gamma$, for every $\gamma \in \Gamma$ (see [Raf05b]). We call $\mathsf{Y}$ the $q$–representative of $Y$.

2.4. **Collar lemma for quadratic differentials.** We recall the following analogue of the collar lemma which relates the quadratic differential lengths of intersecting curves on $S$.

**Theorem 2** ([Raf05b, Theorem 1.3]). *For every $L > 0$, there exists $D_L$, $\log D_L \asymp e^L$, depending only on $L$ such that, if $\alpha$ and $\beta$ are two simple closed curves in $S$ intersecting nontrivially, with $l_\sigma([\beta]) < L$, then*

$$D_L \, l_q([\alpha]) \geq l_q([\beta]).$$

2.5. **Notation.** The notation $x \mathrel{\dot\succ} y$ (respectively, $x \mathrel{\dot\asymp} y$) means that there exists a constant $c$ depending on the topology of $S$ such that $c \, x \geq y$ (respectively, $y/c \leq x \leq c \, y$).

3. SIZE OF A SUBSURFACE

Let $\alpha_Y$ be an essential (non-peripheral) curve in $Y$ with the shortest $q$–length. Define the $q$–size of $Y$ to be $l_q(\alpha_Y)$, and denote it by $\lambda(Y)$, or just $\lambda$ if the subsurface $Y$ is fixed. When $Y$ is a pair of pants (that is, $Y$ has genus 0 and 3 boundary components), there are no essential curves in $Y$. In this case, we define the size of $Y$ to be the maximum $q$–length of boundary components of $Y$. This special case come to play only in the statement of Theorem 4.

We show below that an essential curve in $Y$ whose hyperbolic length is bounded has $q$–length comparable to $\lambda(Y)$. In particular, we are interested in the lengths of curves in a *short marking*, defined as follows.

3.1. **Short markings.** Let $\sigma$ be the hyperbolic metric in the conformal class of $q$, $\epsilon_0$ be a fixed constant smaller than the Margulis constant and $L_0 = \log(1/\epsilon_0)$. Let $\Gamma$ be the set of simple closed $\sigma$–geodesics in $S$ whose $\sigma$–length is less than $\epsilon_0$. An $\epsilon_0$–*thick* component of $\sigma$ is a component $Y$ of $S \setminus \Gamma$. A *short marking* $\mu$ in $Y$ is a set of curves in $Y$ chosen as follows: Take the $\sigma$–shortest pants decomposition of $Y$ (that is a pants decomposition where the sum of $\sigma$-lengths of curves involved is smallest possible) and, for each of these curves, take a transverse curve with the shortest $\sigma$–length. Note that $\mu$ *fills* the subsurface $Y$, that is, all the complementary components of $\mu$ in $Y$ are either disks or annuli sharing a boundary component with $Y$. Also, since the injectivity radius of $Y$ in $\sigma$ is larger than $\epsilon_0$, for every $\beta \in \mu$, $l_\sigma(\beta) \stackrel{.}{\prec} L_0$. However, we think of $\epsilon_0$ as a universal constant and therefore we write $l_\sigma(\beta) = O(1)$.

The marking $\mu$ coarsely determines the geometry of $\sigma$ restricted to $Y$. Let $\mathsf{Y}$ be the $q$-representative of $Y$ and $\lambda = \lambda(Y)$ be the size of $Y$ in $q$. As a first step in comparing the geometry of $(Y, \sigma)$ with $(\mathsf{Y}, |q|)$, we show that $\mu$ is also "short in $\mathsf{Y}$", that is, the $q$–lengths of curves in $\mu$ are comparable with $\lambda$ (the smallest they can be).

**Lemma 3.** *The $q$–lengths of curves in a $\sigma$-short marking $\mu$ of $Y$ are comparable to the size of $Y$. That is, for every $\beta \in \mu$,*
$$l_q([\beta]) \stackrel{.}{\asymp} \lambda.$$

*Proof.* The set $\mu$ fills $Y$; therefore, some curve $\beta \in \mu$ has to intersect $\alpha_Y$. Since $l_q(\alpha_Y) = \lambda$ and $l_\sigma(\beta) = O(1)$ we can conclude from the collar lemma for quadratic differentials (Theorem 2) that $l_q([\beta]) \stackrel{.}{\prec} \lambda$. Similarly, for any curve $\beta'$ that intersects $\beta$, Theorem 2 implies $l_q([\beta']) \stackrel{.}{\prec} \lambda$. But the union of the curves in $\mu$ is a connected set and the number of curves in $\mu$ is bounded by a constant depending on the topology of $S$ only. Therefore, they all have lengths that are bounded by a multiple of $\lambda$, and the multiplicative constant depends on the topology of $S$ only. □

3.2. **The diameter and the area of $\mathsf{Y}$.** We also compare the diameter and area of $\mathsf{Y}$ to $\lambda(Y)$. In general, there is no lower bound for the area of the $q$–representative of a subsurface, as there is in the hyperbolic metric. In fact, there are degenerate cases where the area is zero (see the example at the end of this paper). However, the following statement holds.

**Theorem 4.** *Let $Y$ be a thick component of $\sigma$. Then:*
  (1) $\operatorname{diam}_q(\mathsf{Y}) \stackrel{.}{\asymp} \lambda$.
  (2) $\operatorname{area}_q(\mathsf{Y}) \stackrel{.}{\prec} \lambda^2$.

*Proof.* The $q$–diameter of $\mathsf{Y}$ is larger than half the $q$–length of any essential curve in $\mathsf{Y}$ or any boundary component of $\mathsf{Y}$. Therefore the $q$–diameter is larger than $\lambda/2$. We have to provide an upper bound for the diameter and the area of $\mathsf{Y}$.

Assume $Y$ is not a pair of pants. Consider the marking $\mu$ described in Lemma 3. The $q$–geodesic representatives of curves in $\mu$ divide the surface $\mathsf{Y}$ into disks, punctured disks and annuli (a punctured disk for each puncture in $Y$ and an annulus for each boundary component of $Y$). The numbers of these disks, punctured disks and annuli are bounded by constants depending on the topology of $S$ only. Therefore, it is sufficient to provide an upper bound for the area and the diameter of each component.

Assume $D$ is one of the disk components. Then the $q$-length of $\partial D$ is less than the sum of the $q$-lengths of curves in $\mu$. That is, $l_q(\partial D) \stackrel{.}{\prec} \lambda$. The metric on $D$ is a flat metric with non-positive curvature concentrated on singular points. Therefore, $D$ is a $CAT(0)$ space and satisfies a quadratic isoperimetric inequality (see [BH99, Theorem III.2.17]), that is,

(1) $$\text{area}_q(D) \stackrel{.}{\prec} l_q(\partial D)^2 \stackrel{.}{\prec} \lambda^2.$$

Also, since $D$ is simply connected, non-positive curvature also implies that

$$\text{diam}_q(D) \stackrel{.}{\prec} l_q(\partial D) \stackrel{.}{\prec} \lambda.$$

If $D$ is a component that is a punctured disk, then the total angle around the puncture is only $\pi$ and $D$ does not have singular non-positive curvature. However, if we take a double cover of $D$, the total angle at every point is at least $2\pi$. Now, as in the previous case, the area and the diameter of the double cover of $D$ are bounded above by a multiple of $\lambda^2$ and a multiple of $\lambda$, respectively. Therefore, they also provide upper bounds for the area and the diameter of $D$ as well.

We need to provide an analogous diameter bound and area bound for an annular component $A$. The annulus $A$ has two boundary components: the inner boundary $\gamma_0$, which is a $q$–representative of the core of $A$, and the outer boundary $\gamma_1$. We know that $l_q(\gamma_1) \stackrel{.}{\prec} \lambda$, and, since $\gamma_0$ is shorter than $\gamma_1$, $l_q(\gamma_0) \stackrel{.}{\prec} \lambda$.

The geodesic $\gamma_0$ consists of straight segments which meet at critical points of $q$, where the interior angle at each corner is larger than $\pi$. Let $x_i$, $i = 1, \ldots, n$, be all the corner points where the corresponding interior angles $\theta_i$ are larger than $\pi$. The curvature of $\gamma_0$ with respect to $\mathsf{Y}$, $\kappa_Y(\gamma_0) = \sum_i (\pi - \theta_i)$, is a multiple of $\pi$ (in a quadratic differential, angle modulo $\pi$ is well-defined). Therefore, $\sum_i \theta_i$ is also a multiple of $\pi$. Note that there has to be at least one such corner point because, if $\kappa_Y(\gamma_0) = 0$, one could push $\gamma_0$ into $\mathsf{Y}$ without changing its length. This contradicts the assumption that $\mathsf{Y}$ is disjoint from the interiors of all flat annuli containing geodesic representatives of its boundary components.

Let $A'$ be the double cover of $A$. Denote the lifts of $x_i$ by $y_i$ and $z_i$. Our plan is to fill $A'$ by adding $2(n-1)$ Euclidean triangles to obtain a disk equipped with a singular flat structure that non-positive curvature, that is, in which the total angle at each point is more than $2\pi$, and then use the upper bounds we have for such disks.

We start by attaching a Euclidean triangle to vertices $y_1, y_2, y_3$, which we denote by $\triangle(y_1, y_2, y_3)$ (see Fig. 1). We choose the angle at vertex $y_2$, $\angle y_2$, so that the total angle at $y_2$, $\theta_2 + \angle y_2$, is a multiple of $\pi$. Assuming $0 \leq \angle y_2 < \pi$, there is a unique such triangle. Attach an isometric triangle to $z_1, z_2, z_3$. Now consider points $y_1, y_3, y_4$. Again, there exists a Euclidean triangle with one edge equal to the newly introduced segment $[y_1, y_3]$, another edge equal to the segment $[y_3, y_4]$ and an angle at $y_3$ that makes the total angle at $y_3$, including the contribution from the triangle $\triangle(y_1, y_2, y_3)$, a multiple of $\pi$. Attach this triangle and an identical one to vertices $z_1, z_3, z_4$. Continue in this fashion until finally adding triangles $\triangle(y_1, y_n, z_1)$ and $\triangle(z_1, z_n, y_1)$. Because of the symmetry, the two edges connecting $y_1$ and $z_1$ are equal, and we can glue these together.

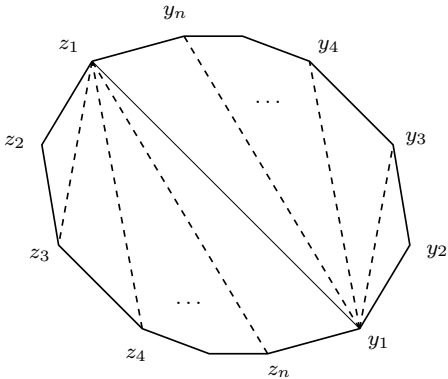

FIGURE 1. **The filling of the annulus $A'$**

For $i \neq 1$, the total angles at $y_i$ and at $z_i$ are multiples of $\pi$ and are larger than $\theta_i > \pi$; therefore, they are at least $2\pi$. We have added $2(n-1)$ triangles. Hence, the sum of the total angles of all vertices is
$$2\sum_i \theta_i + 2(n-1)\pi.$$
which is a multiple of $2\pi$. Therefore, the sum of the angles at $y_1$ and $z_1$ is also a multiple of $2\pi$. But they are equal to each other, and each one is larger than $\pi$. This implies that they are both at least $2\pi$.

We denote this filled annulus by $D'$. The disk $D'$ has a singular flat metric with non-positive curvature, and the $q$–length of the boundary of $D'$ is twice that of $\gamma_1$. The isoperimetric inequality implies that
$$\mathrm{area}_q(D') \prec \lambda^2.$$
This gives the desired upper bound for the area of $A$.

Let $x, y \in A'$, and let $[x, y]$ be the geodesic segment connecting them in $D'$. If $[x, y]$ enters the added part, we can replace the portion of $[x, y]$ that is in the added part with a portion of the inner boundary of $A'$ of length less than or equal to $l_q(\gamma_0) \prec \lambda$. Therefore, the distances in $A'$ are at most $l_q(\gamma_0)$ larger than the distances in $D'$. But $D'$ is simply connected and has non-positive curvature. Therefore, as before, $\mathrm{diam}(D') \prec \lambda$. Hence, $\mathrm{diam}(A') \prec \lambda$. Since $A'$ is a double cover of $A$, the diameter of $A'$ is at most twice that of $A$. This gives the upper bound for the diameter of $A$ that we wanted.

Now assume $Y$ is a pair of pants and let $\alpha$, $\beta$ and $\gamma$ be its boundary components $Y$. Recall that in this case $\lambda(Y) = \max(l_q(\alpha), l_q(\beta), l_q(\gamma))$. The curvature of each boundary component of $Y$ with respect to $Y$ is negative and is a multiple of $\pi$. In particular, it has to be less than $-\pi$. If boundaries of $Y$ do not intersect each other, then the Gauss-Bonnet theorem implies that the sum of the curvatures of boundaries of $Y$ is greater than $2\pi\chi(Y) = -2\pi$, which is a contradiction. Therefore, at least two of the boundaries (say $\beta$ and $\gamma$) have to intersect. We can think of $Y$ as an annulus whose inner boundary is $\alpha$ and whose outer boundary is the curve obtained from concatenation of $\beta$ and $\gamma$. From the above discussion, the $q$–diameter of this annulus is less than a constant multiple of $l_q(\beta) + l_q(\gamma)$ and its $q$–area is less than a constant multiple of $(l_q(\beta) + l_q(\gamma))^2$. But $l_q(\alpha) \leq l_q(\beta) + l_q(\gamma)$ (because $\alpha$ is a $q$–geodesic). Therefore, $\lambda(Y) \mathrel{\dot{\asymp}} l_q(\beta) + l_q(\gamma)$. This finishes the proof in this case. □

## 4. Proof of the main theorem

The idea for the proof of Theorem 1 is to show that the $q$-length of a curve $\alpha$ in a thick component $Y$ can be estimated by the intersection with a short marking $\mu$ in $Y$ as follows

$$(2) \qquad l_q(\alpha) \mathrel{\dot{\asymp}} \lambda(Y) i(\alpha, \mu),$$

and in addition, to use the fact that the hyperbolic length of a curve $\alpha$ in a thick component can be estimated by (see for example [Min93])

$$(3) \qquad l_\sigma(\alpha) \mathrel{\dot{\asymp}} i(\alpha, \mu).$$

To prove Equation(2), we need the following lemma:

**Lemma 5.** *Let $\alpha$ and $\beta$ be two essential curves in $Y$. Then*

$$l_q(\alpha)\, l_q(\beta) \mathrel{\dot{\succ}} \lambda(Y)^2\, i(\alpha, \beta).$$

*Proof.* By perturbing the $q$–geodesic representatives of $\alpha$ and $\beta$, we can find curves $\bar\alpha$ and $\bar\beta$ such that $\bar\alpha$ and $\bar\beta$ intersect at finitely many points (i.e, they do not share a line segment), all the intersections are essential and

$$(4) \qquad l_q(\bar\alpha) \mathrel{\dot{\asymp}} l_q(\alpha) \quad \text{and} \quad l_q(\bar\beta) \mathrel{\dot{\asymp}} l_q(\beta),$$

for constants close to 1. Let $\omega$ be an arc in $\bar\alpha$ of $q$–length $\lambda(Y)/2$ containing the maximum number of intersection points of $\bar\alpha$ and $\bar\beta$. Denote these points by $x_1, ..., x_n$. We have

$$(5) \qquad n \geq \frac{\lambda(Y)/2}{l_q(\bar\alpha)}\, i(\alpha, \beta).$$

These $x_i$ divide $\bar\beta$ into $n$ arcs. The union of any such arc with end points $x_i$ and $x_j$ with the arc in $\omega$ connecting $x_i$ to $x_j$ is an essential curve in $Y$ and therefore has $q$–length greater than $\lambda(Y)$. This implies that the $q$–lengths of these arcs in $\bar\beta$ are all greater than $\lambda(Y)/2$. Therefore,

$$(6) \qquad l_q(\bar\beta) \geq n\, \lambda(Y)/2.$$

Equations (4), (5) and (6) imply the lemma. □

We now prove the main theorem.

**Theorem 6.** *For every essential curve $\alpha$ in $Y$, we have*
$$l_q([\alpha]) \asymp \lambda\, l_\sigma([\alpha]).$$

*Proof.* Let $\mu$ be a short marking on $Y$. By Lemma 5 we have
$$\lambda^2 \sum_{\beta \in \mu} i(\alpha, \beta) \precsim \sum_{\beta \in \mu} l_q(\alpha) l_q(\beta).$$

Applying Lemma 3, we get
$$\lambda^2 \sum_{\beta \in \mu} i(\alpha, \beta) \precsim \#\mu\, \lambda\, l_q(\alpha).$$

But the number of curves $\#\mu$ in $\mu$ depends only on the topology of $S$. Thus, dividing both sides of the above inequality by $\lambda$ and applying Equation(3), we get
$$\lambda\, l_\sigma(\alpha) \precsim l_q(\alpha).$$

To prove the inequality in the other direction, recall from the proof of Theorem 4 that the complementary components of $q$–geodesic representatives of curves in $\mu$ in $Y$ have diameter bounded above by a multiple of $\lambda$. Therefore, we can retrace the intersection pattern of $\alpha$ with $\mu$ and produce a curve in $Y$ in the homotopy class of $\alpha$ whose length is bounded by a multiple of $\lambda \sum_{\beta \in \mu} i(\alpha, \beta)$. Therefore,
$$\lambda\, l_\sigma(\alpha) \succsim l_q(\alpha).$$

This completes the proof. □

## 5. An example

In this final section, we construct an example of a subsurface whose $q$–area is zero and calculate its diameter.

Consider the singular Euclidean metric on a genus 4 surface obtained by cutting 4 parallel slits on a flat torus of area 1 and gluing the boundaries as in Fig. 2.

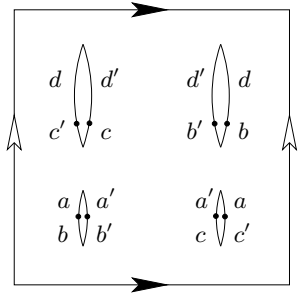

Figure 2. **Singular Euclidean metric on $S$**

Topologically, $S$ can be obtained by gluing a 4-times punctured sphere $Y$ to a 4–times punctured torus (the flat torus we started with). The Euclidean metric on the flat torus defines a conformal structure on $S$ and, since the slits are parallel, choosing an angle for the leaves of horizontal foliation defines a compatible quadratic differential $q$ on $S$. The $q$-representative $\mathsf{Y}$ of $Y$ in this case is degenerate and has area zero. In fact, it is a graph with 2 vertices and 4 edges whose Euclidean lengths are equal to the Euclidean lengths of segments $a$, $b$, $c$ and $d$ (Fig. 3, left).

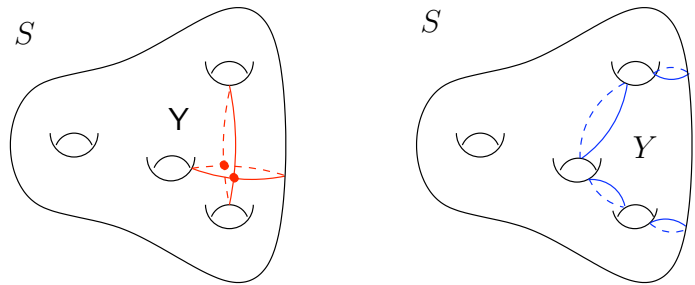

FIGURE 3. **Subsurface $Y$ and its $q$–representative $\mathsf{Y}$**

Let $l_1$ be the $q$–length of the top two slits and $l_2$ be the $q$–length of the bottom two slits ($l_1 > l_2$). By choosing $l_1$ and $l_2$ small enough and far apart, we can ensure that the hyperbolic lengths of boundaries of $Y$ are smaller than $\epsilon_0$. (In fact, this would allow a round annulus of large modulus to be embedded around each boundary curve.) Note that $\diam_q(\mathsf{Y}) = l_1$ and $\lambda(Y) = 2\, l_2$. In this case, Theorem 4 implies that, if $Y$ is a thick piece of $S$, that is, if no essential curves in $Y$ have hyperbolic lengths less than $\epsilon_0$, then $l_1$ and $l_2$ are almost equal. Curves in $Y$ can be parametrized with rational numbers. Let $\alpha$ be the $\frac{m}{n}$ curve in $Y$. Then, $\alpha$ has a $q$–length of $2m\, l_1 + 2n\, l_2$ and a $\sigma$–length comparable to $m + n$. If $l_1 \asymp l_2$, we see (as predicted by Theorem 6) that
$$l_q([\alpha]) = 2m\, l_1 + 2n\, l_2 \asymp l_2(m+n) \asymp \lambda(Y)\, l_\sigma([\alpha]).$$